\documentclass[a4]{article}
\usepackage{amsmath}
\usepackage{amsthm}
\usepackage{eufrak}

\begin{document}
\title{On statistics of permutations  chosen from  the Ewens distribution}

\author{Tatjana Bak\v sajeva\thanks{ Vilnius University, Naugarduko str. 24, LT-03225 Vilnius,
Lithuania} and Eugenijus Manstavi\v cius\thanks{Vilnius University, Institute and Faculty of Mathematics and Informatics,
Akademijos 4, LT-08663, Vilnius, Lithuania} }

\maketitle

\footnotetext{{\it AMS} 2000 {\it subject classification.} Primary
60C05;      secondary 05A16, 20P05. \break {\it Key words and
phrases}. Random permutation, cycle structure,
  additive function, multiplicative function, permutation matrix, Poisson law, quasi-Poisson law.}

%\address{Department of Mathematics and Informatics, Vilnius University, Naugarduko str. 24, LT-03225 Vilnius,  Lithuania}

\begin{abstract}
We explore the asymptotic distributions of sequences
of integer-valued  additive functions defined on the symmetric group endowed with
the Ewens probability measure as the order of the group increases.
 Applying the method of factorial moments, we establish necessary and sufficient conditions for the weak convergence of distributions to discrete laws.
 More attention is paid to the Poisson limit distribution. The particular case of the number-of-cycles with restricted lengths function is analyzed in
 more detail. The results can be applied to statistics defined on random
permutation matrices.
\end{abstract}

%\maketitle

\newtheorem{thm}{Theorem}
\newtheorem{lem}{Lemma}
\newtheorem{cor}{Corollary}
\newtheorem*{cor*}{Corollary}
\newtheorem{conj}{Conjecture}

\newtheorem{prop}{Proposition}
\newtheorem*{prop*}{Proposition}
\newtheorem*{CLT}{CLT (\cite{EM-LMJ96})}
\newtheorem*{Comp}{Compactness Thm (\cite{EM-RJ08})}

\newtheorem*{WLLN}{WLLN (\cite{EM-RJ08})}

\def\E{\mathbf{E}}
\def\C{\mathbf{C}}
\def\D{\mathbf{D}}
\def\G{\mathbf{G}}
\def\V{\mathbf{V}}
\def\Ra{\Rightarrow}
\def\N{\mathbf{N}}
\def\R{\mathbf{R}}
\def\S{{\mathbf{S}_n}}
\def\Z{\mathbf{Z}}

\def\k{\kappa}

\def\e{\varepsilon}
\def\ro{{\rm o}}

\def\rO{{\rm O}}
\def\re{{\rm e}}
\def\rd{{\rm d}}

\def\sj{\sum_{j\le n}}
\def\n{$n\to\infty$}
\def\cF{\mathcal F}
\def\s{\smallskip}

\section{Introduction}

   We deal with asymptotic value distribution problems of mappings defined on the symmetric group $\S$ as \n.
    Let $\sigma\in\S$ be an arbitrary permutation and
   $
    \sigma=\k_1\cdots\k_w
   $
be its representation as the product of independent cycles $\k_i$
and $w:=w(\sigma) $  be their number. If $k_j(\sigma)$, $1\leq
j\leq n$, denotes the number of cycles of length $j$ in this
decomposition, then $\bar k(\sigma):=\big( k_1(\sigma), \dots,
k_n(\sigma)\big)$ is called the cycle structure vector. The Ewens
probability measure on the subsets $A\subset \S$ is defined by
   \[
     \nu_n(A):=\nu_{n,\theta}(A)={1\over \theta^{(n)}} \sum_{\sigma\in A}\theta^{w(\sigma)},
     \]
     where $\theta>0$ is a fixed parameter and $\theta^{(n)}:=\theta(\theta+1)\cdots(\theta+n-1)$.
An easy combinatorial argument gives the  distribution of the cycle structure vector
      \begin{equation}
  \nu_n\big(\bar k(\sigma)=\bar
  s\big)={n!\over \theta^{(n)}}\prod_{j=1}^n \Big({\theta\over j}\Big)^{s_j} {1\over s_j!},
  \label{nub}
  \end{equation}
  where
 $\bar s=(s_1,\dots, s_n)\in\Z_+^n$ and $\ell(\bar s):=1s_1+\cdots+ns_n=n$. The probability (\ref{nub}), ascribed to $\bar s\in \ell^{-1}(n)\subset\Z_+^n$, i.e.
 \[
     P_\theta(\{\bar s\})={n!\over \theta^{(n)}}\prod_{j=1}^n \Big({\theta\over j}\Big)^{s_j} {1\over s_j!}
 \]
 is  called the Ewens Sampling Formula. It defines a probability measure  on $\ell^{-1}(n)$  (afterwards, we denote it by ESF).
 If  $\xi_j$, $j\geq 1$, denote independent Poisson r.vs  given on some
probability space $\{\Omega,\mathcal F, P\}$ with
$\E\xi_j=\theta/j$ and $\bar \xi:=(\xi_1,\dots, \xi_n)$, then
\begin{equation}
  \nu_n\big(\bar k(\sigma)=\bar s\big)=  P\big(\bar \xi=\bar s \big|\,\ell(\bar \xi)=n\big), \quad \bar s\in\ell^{-1}(n).
  \label{Cond}
  \end{equation}
  Moreover, the total variation distance
  \begin{equation}
  {1\over 2}\sum_{s_1,\dots, s_r\geq0}\big|\nu_n\big(k_1(\sigma)=s_1,\dots,
  k_r(\sigma)=s_r\big)-P(\xi_1=s_1,\dots,
  \xi_r=s_r\big)\big|=o(1)
  \label{totalvar}
  \end{equation}
if and only  if $r=o(n)$. Here and in what follows we assume that \n.
The book  \cite{ABT} is a good reference for the listed and many more properties of the ESF.

In the present paper, we discuss the asymptotic value distribution of an additive  (completely additive) function $h:\S\to\R$ with respect to $\nu_n$.
 By definition, such a function  is  defined via a real array  $\big\{a_{j}, j\geq 1\big\}$,  by setting
\begin{equation}
   h(\sigma):=\sum_{j\leq n} a_{j}k_j(\sigma).
 \label{h}
 \end{equation}
By virtue of $\ell\big(\bar k(\sigma)\big)=n $ if $\sigma\in \S$,
$h(\sigma)$, as a r.v.  under the probability measure $\nu_n$, is
a sum of dependent summands.  We prefer to leave the elementary
event $\sigma$ in its notation. This goes in some contrast to
other r.vs defined as above on a nonspecialized space
$\{\Omega,\mathcal F, P\}$.  Taking arrays $a_{nj}$, $1\leq j\leq
n$, $n\geq 1$, we obtain sequences of functions $h_n(\sigma)$. So,
if $a_{nj}=1$ for $j\in J_n\subset \{1,\dots, n\}$ and $a_{nj}=0$
elsewhere, we have a sequence of the number-of-cycles with
restricted lengths additive functions which we denote by
$w(\sigma,J_n)$.

 Apart from the latter instance, the additive functions are involved in many combinatorial, algebraic and statistical problems.
 The function $h(\sigma)$ defined via $a_j=\log j$, $j\leq n$, well approximates the logarithm of group theoretical order of almost all
  permutations $\sigma\in \S$ (see \cite{ET-ZWar65}  and \cite{VZ-LMJ02} or \cite{VZ-LMJ04}).
  Particular additive functions appear in physical models as a part of Hamiltonians in the Bose gas  theory (see \cite{BeUel-CMPh09}--\cite{BeUel-EJP11}
  and the references therein).
  They are indispensable treating the  random permutation matrix ensemble.
 Let $M:=M(\sigma):= \big({\mathbf 1}\{i=\sigma(j)\}\big)$, $1\leq i,j\leq n$ and $\sigma\in\S$, be such a matrix taken with the weighted frequency $\nu_n(\{M\})=\nu_n(\{\sigma\})=\theta^{w(\sigma)}/\theta^{(n)}$,
      \begin{equation}
   Z_n(x;\sigma):=\operatorname{det}\big(I-xM(\sigma)\big)=\prod_{j\leq n}(1-x^j)^{k_j(\sigma)}
   \label{Tr}
   \end{equation}
   be its characteristic polynomial, and let $\re^{2\pi i \varphi_j(\sigma)}$, where  $\varphi_j(\sigma)\in [0,1)$ and $j\leq n$,
   be its eigenvalues. The papers \cite{HKCS-SPA00}, \cite{KW-AP00}, \cite{KW-JTP03}, and \cite{DZ-EJP10}  or  some preprints put in the AMS arXiv (see, for instance,  \cite{BArKDa-Ar11} and \cite{HNNZ-ar11}  and the references therein) concern  $\log|Z_n(x;\sigma)|$, $\Im \log Z_n(x;\sigma)$
   or the trace related statistics
   \begin{equation}
      \operatorname{Trf}(\sigma):=\sum_{j\leq n} f\big(\varphi_j(\sigma)\big)=\sum_{j\leq n} k_j(\sigma)
      \sum_{0\leq s\leq j-1}f\Big({s\over j}\Big),
      \label{Trf}
      \end{equation}
      where $f:[0,1]\to \R$ is a  function. An indicator function $f={\mathbf 1}_A$ of an interval $A\subset [0,1]$ or other integer-valued functions fall within our objectives.

So far, the general problem to find necessary and sufficient conditions under which the distribution function
\[
         T_n(x):=\nu_n\big(h_n(\sigma)<x\big)
         \]
weakly converges to a limit law is out of reach even for
$\theta=1$. The case when the sequence $a_{nj}=a_j$, $j\geq 1$,
that is, does not depend on $n$, is easy. Then the answer is given
by an analog of the three series theorem of Kolmogorov (see
Theorem 8.25 in \cite{ABT}). If $\theta=1$ and
$h_n(\sigma)=h(\sigma)/\beta(n)$, where the function $h(\sigma)$
is fixed and $\beta(n)>0$, $\beta(n)\to\infty$ but
$\beta(un)/\beta(n)\to 1$ for every fixed  $0<u<1$ (slowly
oscillating at infinity), necessary and sufficient conditions were
established in the second author's paper \cite{EM-LMJ11} which
contains an extensive reference list of earlier papers by other
authors. If $\beta(n)$ is regularly varying at infinity, the first
results go back to  paper \cite{EM-LMJ96}. The possed problem
remains  open if no \textit{a fortiori} condition on $\beta(n)$ is
taken. On the other hand, for partial sum processes defined by
additive functions, convergence of distributions in appropriate
function spaces to infinitely divisible measures implies slow
oscillation of  $\beta(n)$. This further yields necessary and
sufficient convergence conditions even for generalized Ewens
probability measures. For the latest account in this direction, we
refer to  \cite{BoEM-Ku12}.

 In the present paper, we focus on  sequences of additive functions $h_n(\sigma)$ defined via $a_{nj}$. Afterwards, we will use the abbreviation $a_j=a_{nj}$ without the index $n$ and take $a_j=0$ if $j>n$. If $\theta=1$, the partial sums of such functions have been used to model stochastic processes \cite{EM-APM02}. Recently, we \cite{TKEM-AUBud13} succeeded  to establish necessary and sufficient conditions for the weak law of large numbers if $\theta\geq 1$. Some success in proving general limit theorems has been  achieved for the integer-valued functions $h_n(\sigma)$. If $\theta=1$, this case has been explored by the second author in \cite{EM-LMJ05}--\cite{EM-RJ08}. For $\theta>0$, the first author in  \cite{TK-SMS07} and \cite{TK-LMD09} obtained an exhaustive result for the sequence $w(\sigma,J_n)$. We now generalize this dealing with the case $a_j\in \Z_+$ if $j\leq n$ and supply a few instances shedding more light about the class of possible limit distribution for $w(\sigma,J_n)$. On the other hand, one of the purposes of the present paper is to demonstrate the factorial moment method.
  The approach proved to be useful in a series of the number-theoretical papers by J. \v Siaulys \cite{JS-LMJ96}-\cite{JS-LMJ00}. The idea lays in analysis of the  expressions of moments. Though involved, they contain the key information useful in establishing necessary and sufficient conditions for the convergence of distributions.

In what follows, let  $\Rightarrow$ stand for weak convergence and
$F_Y(x)$ be a distribution function of a r.v. $Y$ concentrated on
$\Z_+=\N\cup\{0\}$. The mean value with respect to $\nu_n$ of a
function $g(\sigma)$ defined in $\S$ will be denoted by $\E_n
g(\sigma)$. Set $x_{(r)}=x(x-1)\cdots (x-r+1)$, $r\geq 1$ for the
falling factorial and $x_{(0)}=1$. Let $\rho_n(m)$ be an error
term, not the same in different places but satisfying the relation
\[
   \lim_{m\to \infty} \limsup_{n\to \infty} |\rho_n(m)|=0.
\]
Afterwards, if this is not indicated, we take $i, r, r_i, j, j_i\in \N$ and $j, j_i\leq n$.  The first two theorems involve the quantity
\begin{eqnarray*}
\Upsilon_n(l,m)&:=& \sum_{u=1}^l\theta^u\sum_{r_1+\cdots+r_u=l\atop 1\leq r_i\leq m, i\leq u}\bigg({{l-1}\atop{r_1-1}}\bigg)\cdots
      \bigg( {{l-r_1-\cdots-r_{u-1}-1}\atop{r_u-1}}\bigg)\nonumber
      \\&&\times\sum_{j_1+\cdots+j_u< n}{ a_{j_1(r_1)}\cdots
       a_{j_u(r_u)}\over{j_1\cdots j_u}}\Big(1 - \frac{j_1 + \cdots + j_u}{n}\Big)^{\theta - 1}
\end{eqnarray*}
which is an approximation of the $l$th factorial moment of an appropriately  truncated additive function obtained from $h(\sigma)$. The first result concerns necessary conditions for convergence.

\begin{thm}
\label{thm-but} Let $\theta\geq 1$, $h_n(\sigma)$ be a sequence of integer-valued additive functions, and let $Y$ be a r.v. taking  values in $\Z_+$ and such that $\E Y^\alpha<\infty$ for  $\alpha\geq 2+\e>2$.
 If  $T_{n}(x)\Rightarrow  F_Y(x)$, then
 \begin{equation}
    \sum_{j<n} {{\mathbf 1}\{a_j\leq -1\}\over j} \Big(1-{j\over n}\Big)^{\theta-1}=\ro(1)
    \label{ajneg}
    \end{equation}
    and
  \begin{equation}
                      \Upsilon_n(l,m)-\E Y_{(l)}=\rho_n(m)
\label{Yps}
 \end{equation}
for each fixed natural number $l \leq \alpha -1-\e $.
\end{thm}

We have to confess  that the technical condition $\theta\geq 1$, in Theorem \ref{thm-but} and in some subsequent results concerning the necessity, is undesirable. Sufficient convergence conditions are   given by the following result.

 \begin{thm} \label{thm-pak} Let  $\theta>0$ and $h_n(\sigma)$ be a sequence of integer-valued additive functions. Assume that condition $(\ref{ajneg})$
 is satisfied. If there exists a sequence $\Upsilon(l)$ such that
 \begin{equation}
                      \Upsilon_n(l,m)-\Upsilon(l)=\rho_n(m)
\label{Yps1}
 \end{equation}
 for every $l\in\N$ and
 \[
      \sum_{l=0}^\infty {\Upsilon(l) 2^l\over l!}<\infty,
 \]
 then  $T_n(x)\Rightarrow  F_Y(x)$ and  $\E Y_{(l)}=\Upsilon(l)$ for $l\geq 1$.
 \end{thm}

The following  corollary gives necessary and sufficient convergence conditions for more specialized cases. Let $\Pi_\mu(x)$ be the distribution function of the Poisson law with parameter $\mu>0$.

\begin{cor} \label{Poisson} Let $\theta\geq1$, $h_n(\sigma)$ be a sequence of integer-valued additive functions. The convergence $T_{n}(x)\Rightarrow \Pi_\mu(x)$ holds if and only if condition $(\ref{ajneg})$ holds and
\begin{equation}
                              \Upsilon_n(l,m)-\mu^l=\rho_n(m)
\label{Ypsmu}
 \end{equation}
 for every $l\in\N$.
\end{cor}

As the total variation estimate approximation (\ref{totalvar})
shows that the Poisson distribution appears as a limit if the
cycles of lengths up to $r=o(n)$ are involved. By the next
corollary, we demonstrate that one can find  $a_j=a_{nj}$, $n/2<
j\leq n$, which defines a sequence of additive functions obeying
a Poisson limit law with a sufficiently small $\mu$. Such a phenomenon has been observed in \cite{EM-APM02} if $\theta=1$.
  The construction involves the  following strictly increasing in $x\in [1/2,1]$ function
\[
   t_\theta(x):=\theta\int_{1/2}^x (1-u)^{\theta-1} { du\over u}.
\]
We will prove that  $t_\theta(1)<1$ if $\theta\geq1$.

 \begin{cor} \label{cor-n-puse} Let $\theta\geq1$, $\mu\leq -\log(1-t_\theta(1))$, and $a_j\in \Z_+$ so that
\begin{equation}
   \sum_{j\leq n/2} {{\mathbf 1}\{a_j\not=0\}\over j}=o(1).
  \label{nupj2}
  \end{equation}
    The convergence $T_n(x)\Rightarrow \Pi_\mu(x)$ holds if and only if
\begin{equation}
\theta\sum_{n/2<j<n} {{\mathbf 1}\{a_j=k\}\over j}\Big(1-{j\over n}\Big)^{\theta-1}= {\rm e}^{-\mu}{\mu^k\over k!}+o(1)
\label{poisk}
\end{equation}
for every fixed $k\in\N$.
    \end{cor}

 An instance satisfying (\ref{poisk})  will be provided below after the proof of this corollary.

 The next part of the paper deals with additive functions defined via bounded $a_j$ for the overwhelming proportion of $j\leq n$. This allows us to obtain some results for $\theta<1$.

\begin{thm} \label{thm-aprez} Let $\theta>0$, $a_j\in\Z_+$, $j\leq n$, and, for some $K\in \N$,
\begin{equation}
 \sum_{j <n} {{\mathbf 1}\{a_j\geq K\}\over j}\Big(1-{j\over n}\Big)^{\theta-1}= o(1).
 \label{aprez}
 \end{equation}
 The convergence $T_{n}(x)\Rightarrow  F_Y(x)$ holds if and only if there exist a sequence $\Upsilon(l)$ such that
\begin{equation}
\lim_{n\to\infty}\Upsilon_n(l)=\Upsilon(l).
\label{Up-apr}
 \end{equation}
 for every $l\in\N$. If this condition holds, $\E Y_{(l)}=\Upsilon(l)$ for $l\geq 1$.
\end{thm}

For the number-of-cycles with restricted lengths functions $w(\sigma,J_n)$, condition (\ref{Up-apr}) attains the most simple form.
Let the asterisk $\ast$ over a  sum substitute for  the condition $a_j=1$ or, equivalently, for $j\in J_n$.
 Set
\[
V_n(x):= \nu_n\big( w(\sigma, J_n)<x\big)
\]
and
\[
\upsilon_n(l):=\theta^l\sum_{j_1+\cdots+j_l<n}^\ast{1\over{j_1\cdots j_l}}\Big(1 - \frac{j_1 + \cdots + j_l}{n}\Big)^{\theta - 1}, \quad l\geq 1.
\]

The next corollary of Theorem \ref{thm-aprez} has been proved in \cite{TK-SMS07} and \cite{TK-LMD09}.

\begin{cor} \label{cor-Y} Let $\theta>0$ and $J_n\subset \{1,\dots, n\}$ be arbitrary. The distribution function
$V_n(x)  \Rightarrow  F_Y(x)$
  if and only if there exists a sequence $\upsilon(l)$, $l\geq 1$, such that
  \begin{equation}
\lim_{n\to\infty}\upsilon_n(l)=\upsilon(l)
\label{upsilon}
 \end{equation}
 for every $l\in\N$. If the latter condition is satisfied, then $\E Y_{(l)}=\upsilon(l)$ for $l\geq1$.
\end{cor}

Here is a particular case.

\begin{cor} \label{cor-Poi} Let $\theta\geq 1$ and $J_n\subset \{1,\dots, n\}$ be arbitrary. The convergence $V_{n}(x)\Rightarrow  \Pi_\mu(x)$ holds if and only if there exists a sequence $r=r(n)=o(n)$ such that condition
\begin{equation}
\sum_{j\leq r}^\ast {\theta\over j}=\mu + o(1),\qquad
\sum_{r<j< n}^\ast {1\over j}\Big(1-{j\over n}\Big)^{\theta-1}= o(1).
\label{Poi-tr}
\end{equation}
 is satisfied.
\end{cor}

Corollary  \ref{cor-Poi} demonstrates that, for $w(\sigma,J_n)$,  the limit Poisson law can be supported by short cycles only. If $\theta=1$, this fact has been observed in \cite{EM-RJ08} even  for bounded  $a_j$ where $j\leq n$.
Non-degenerate limit distributions concentrated on the finite set $\{0,1,\dots, L-1\}$ where $L\geq 2$ raise a particular interest.

\begin{cor} \label{cor-finsup}
Let $\theta>0$, $J_n\subset \{1,\dots, n\}$, and $L\in\N\setminus\{1\}$ be arbitrary. Assume that  $Y$ is a r.v. taking values in
$\{0,1,\dots, L-1\}$ and such that  $\E Y_{(l)}=\upsilon(l)$ if $l\leq L-1$. The convergence
$V_{n}(x)\Rightarrow F_Y(x)$ holds under the  following conditions:
\begin{equation}
\lim_{n\to\infty}\theta^l\sum_{n/L<j_1,\dots,j_l<n}^\ast{{\mathbf 1}\{j_1+\cdots+j_l<n\}\over{j_1\cdots j_l}}\Big(1 - \frac{j_1 + \cdots + j_l}{n}\Big)^{\theta - 1}=\upsilon(l)
\label{up-nL}
\end{equation}
 for each $l\leq L-1$ and
 \begin{equation}
             \sum_{j\leq n/L}^\ast{1\over j}=o(1).
 \label{tr-nL}
 \end{equation}

 Conversely, if $\theta\leq 1$ and $V_{n}(x)\Rightarrow F_Y(x)$, then conditions $(\ref{up-nL})$ and  $(\ref{tr-nL})$ are satisfied.
 \end{cor}

M. Lugo \cite{ML-ArX09} discussed  a case  falling within the scope of the last corollary. By definition,
a r.v. $\mathfrak q$ is called $(k,\mu)$ quasi-Poisson if  it  has the distribution
\[
                   P({\mathfrak q}=i)=\sum_{j=i}^{k} {j\choose i} (-1)^{j-1} \lambda^j, \quad i=0,1,\dots k,
\]
where $0<\lambda\leq 1$. The factorial moment $\E{\mathfrak
q}_{(l)}=\lambda^l$ if $l\leq k$ and $\E {\mathfrak q}_{(l)}=0$ if
$l>k$.
If $\theta=1$, it is easy to define a subset $J_n$ so that $w(\sigma,J_n)$ obeys the quasi-Poissonian limit law.
 Actually, this and some other results from \cite{ML-EJC09} were already  contained in Theorem 1.3 of the second author's paper
\cite{EM-AMUO05}. Lugo wrote on page 13 of \cite{ML-ArX09}: \textit{...in the case of the Ewens distribution, the
following conjecture seems reasonable}:

\s

{\bf Conjecture 15.} \textit{ The expected number of cycles of length in} $[\gamma n, \delta n]$ \textit{ of a permutation of} $\{1,\dots, n\}$ \textit{chosen from the Ewens distribution approaches}
\[
   \lambda= \int_\gamma^\delta (1-x)^{\theta-1} {dx\over x}
\]
\textit{as} $n\to\infty$. \textit{Furthermore, in the case where $1/(k+1)\leq \gamma<\delta<1/k$
for some positive integer $k$, the distribution of the number of cycles
converges in distribution to quasi-Poisson $(k,\lambda)$.}

\s

As we will see, the factor $\theta$ is missing in the formula for $\lambda$ and, if $\theta\not=1$, the claim of Hypothesis is false. The limit law for
 Lugo's instance  does exist but is not quasi-Poisson. To see this, it suffices to approximate the moments $\upsilon_n(l)$ by appropriate $l$-fold integrals
 and to check that some of the relations $\upsilon(l)=\upsilon(1)^l$, where $l\leq k$, fails. The details are given in the last section.

The paper is organized as follows. Section 2 contains formulae of
factorial moments and needed lower estimates of some frequencies.
The proofs of Theorems \ref{thm-but} and \ref{thm-pak} are
presented in Sections 3 and 4. The next section deals with the
case of bounded $a_j$ including also $w(\sigma,J_n)$. The last
Section gives some illustrative instances. Some distributions
which never appear as limits for $V_n(x)$ are also indicated.

\section{Lemmata}

In this section, we present exact expressions of the factorial moments of a completely additive functions $h(\sigma)$ defined via  $a_j \in \R$.  Particular attention is spared to the case with bounded $a_j $ and  approximations. Denote
\[
\psi_n(m)={n!\over \theta^{(n)}}\, {\theta^{(m)}\over
 m!}= \prod_{k=m+1}^n \Big(1+{\theta-1\over k}\Big)^{-1},
 \]
  where $0\leq m\leq n$. It is well known that
 \[
    {\theta^{(m)}\over m!}= {m^{\theta-1}\over \Gamma(\theta)}\bigg(1+O\Big({1\over m}\Big)\bigg), \quad m\geq 1,
    \]
    where $\Gamma(u)$ is Euler's gamma-function. Hence
    \begin{equation}
     \psi_n(m)=\Big({m\over n}\Big)^{\theta-1}\bigg(1+O\Big({1\over m}\Big)\bigg),\quad  1\leq m\leq n.
\label{psinm}
\end{equation}
In the sequel, we will use the  inequalities
\[
      \psi_n(n-i-j)\geq  \psi_n(n-i) \psi_n(n-j) \quad {\rm if}\quad  \theta\leq 1
\]
and
\[
      \psi_n(n-i-j)\leq  \psi_n(n-i) \psi_n(n-j) \quad {\rm if}\quad \theta\geq 1,
\]
valid for $0\leq i,j\leq n$.

 It is worth to recall  Watterson's formula.

  \begin{lem}
\label{lem-Wat} For $(j_1,\cdots,j_r) \in \mathbf{Z}_{+}^r$, $l = 1j_1 + \cdots + rj_r$ and $1 \leq r \leq n$,
\begin{equation}
\E_{n}\big\{\prod_{i=1}^r {k_i}_{(j_i)}(\sigma)\big\} =\psi_n(n-l)\mathbf{1}\{l\leq n\}
\prod_{i=1}^r\bigg(\frac{\theta}{i}\bigg)^{j_i}.
\label{Wat}
\end{equation}
\end{lem}

{\it Proof}. See  (5.6) on page 96 in \cite{ABT}.

\s

 The next lemma extends in some way the previous formula.

\begin{lem}
\label{FM}
Let $\theta>0$. For a completely additive function $h(\sigma)$   and every $k\in\N$, we have
\begin{eqnarray}
&&\mathbf{E}_{n} h(\sigma)_{(k)}=\gamma_n(k)\nonumber\\
&:=& \sum_{u=1}^k\theta^u\sum_{r_1+\cdots+r_u=k}{k-1\choose r_1-1}\cdots
      {k-r_1-\cdots-r_{u-1}-1 \choose r_u-1}\nonumber\\
      &&
      \quad \times\sum_{j_1+\cdots+j_u\leq n}{ a_{j_1(r_1)}\cdots
       a_{j_u(r_u)}\over{j_1\cdots j_u}}\psi_n\big(n-(j_1 + \cdots + j_u)\big).
\label{gammank}
\end{eqnarray}
\end{lem}

 {\it Proof.} We first prove a recurrence relation for $\beta_n(k):= \big(\theta^{(n)}/{n!}\big)\gamma_n(k)$. Set
$\beta_n(0)=\theta^{(n)}/ n!$ if $n\geq 0$. Moreover, let $\beta_0(k)=0$ if $k\geq 1$.
 Further, let
$\varphi_0(z)=1$  and
  \[
  \varphi_n(z)={\theta^{(n)}\over n!}\mathbf{E}_{n} z^{h(\sigma)}.
\]
 Thus, $\varphi_n^{(k)}(z)|_{z=1}=\beta_n(k)$.

 Grouping over the classes of $\sigma$ with the common cycle vector and using Cauchy's formula for the cardinality of a class, we have
    $$
    \varphi_n(z)=\sum_{\sigma\in\mathbf{S}_n}
\prod_{j=1}^n\bigg({\theta z^{a_j}\over j}\bigg)^{k_j}{1\over{k_j!}}.
  $$
  This leads to the formal series equality
  $$
  \sum_{n\geq 0}\varphi_n(z)w^n=\exp\bigg\{\theta\sum_{j\geq1}{z^{a_j}\over j} w^j\bigg\}
  $$
  and
  \begin{eqnarray*}
  \sum_{n\geq0}\varphi_n'(z)w^n&=&
  \theta\sum_{m\geq 0}\varphi_m(z)w^m\cdot \sum_{j\geq1}{a_jz^{a_j-1} \over j} w^j\\
  &=&
  \theta\sum_{n\geq 0}\bigg(\sum_{j\leq n}\varphi_{n-j}(z){a_j z^{a_j-1}\over j}\bigg) w^n.
  \end{eqnarray*}
Hence
\begin{equation}
\varphi_n'(z)=
  \theta\sum_{j\leq n}\varphi_{n-j}(z){a_j z^{a_j-1}\over j}.
\label{deriv}
\end{equation}
Taking the derivatives with respect to $z$ of the $(k-1)$th order, we arrive at
    $$
   \varphi_n^{(k)}(z)=\theta\sum_{j\leq n}\sum_{l=0}^{k-1}\bigg({{k-1}\atop{l}}\bigg){a_{j(l+1)}z^{a_j-l-1}\over j}
   \varphi_{n-j}^{(k-1-l)}(z).
$$
Consequently, if $z=1$, we obtain
\begin{equation}
   \beta_n(k)=\theta\sum_{r=1}^{k-1}\bigg({{k-1}\atop{r-1}}\bigg)\sum_{j\leq n}{a_{j(r)}\over j}
   \beta_{n-j}(k-r)+\theta\sum_{j\leq n}{a_{j(k)}\over j}{\theta^{(n-j)}\over (n-j)!}.
\label{recurr}
\end{equation}

We now apply the mathematical induction to prove that
\begin{eqnarray}
\beta_n(k)&=& \sum_{u=1}^k\theta^u\sum_{r_1+\cdots+r_u=k}{k-1\choose r_1-1}\cdots
      {k-r_1-\cdots-r_{u-1}-1 \choose r_u-1}\nonumber \\
      &&
      \quad \times\sum_{j_1+\cdots+j_u\leq n}{ a_{j_1(r_1)}\cdots
       a_{j_u(r_u)}\over{j_1\cdots j_u}}{\theta^{(n-j_1- \cdots - j_u)} \over (n-j_1- \cdots - j_u)!}.
\label{rec-beta}
\end{eqnarray}
A direct application of (\ref{deriv}) yields
\[
    \beta_n(1)=\theta\sum_{j\leq n} {a_j\over j} {\theta^{(n-j)}\over (n-j)!}.
    \]

Assume that the induction hypothesis (\ref{rec-beta}) holds for
$\beta_{n-j}(k-r)$  if $k-r\geq1$.   Applying this  formula, we use  the summation indexes $r_2,\dots$ and $j_2,\dots$ leaving $r_1$ and $j_1$ for the summation in  (\ref{recurr}) with respect to $r$ and $j$. So, inserting the assumption into (\ref{recurr}), we obtain
\begin{eqnarray*}
  &&\beta_{n}(k)=\theta\sum_{r_1=1}^{k-1}
  \bigg({{k-1}\atop{r_1-1}}\bigg)
  \sum_{j_1\leq n}{a_{j_1(r_1)}\over j_1}\\
  &&\times
   \sum_{u=2}^{k-r_1+1}\theta^{u-1}\sum_{r_2+\cdots+r_u=k-r_1}
  \bigg({{k-r_1-1}\atop{r_2-1}}\bigg)\cdots
      \bigg( {{k-r_1-\cdots-r_{u-1}-1}\atop{r_u-1}}\bigg)\\
      &&
      \times\sum_{j_2+\cdots+j_u\leq n-j_1}{a_{j_2(r_2)}\cdots
       a_{j(r_u)}\over{j_2\cdots j_u}}\frac{\theta^{(n-j_1-\cdots-j_u)}}{(n-j_1-\cdots-j_u)!}
       +\theta\sum_{j=1}^n{a_{j(k)}\over{j}}\frac{\theta^{(n-j)}}{(n-j)!}.
       \end{eqnarray*}
  Interchanging the summation, we arrive at
  \begin{eqnarray*}
      &&\beta_n(k)      =
        \sum_{u=2}^k\theta^u\sum_{r_1+\cdots+r_u=k}\bigg({{k-1}\atop{r_1-1}}\bigg)\cdots
      \bigg( {{k-r_1-\cdots-r_{u-1}-1}\atop{r_u-1}}\bigg)\\
      &&
      \times\sum_{j_1+\cdots+j_u\leq n}{a_{j_1(r_1)}\cdots a_{j_u(r_u)}\over{j_1\cdots j_u}}\frac{\theta^{(n-j_1-\cdots-j_u)}}{(n-j_1-\cdots-j_u)!}+
       \theta\sum_{j=1}^n{a_{j(k)}\over{j}}\frac{\theta^{(n-j)}}{(n-j)!}.
\end{eqnarray*}

The last sum equals the summand corresponding to $u=1$  in the previous sum over $u$. Joining them together we obtain (\ref{rec-beta}). Further, dividing it by $\theta^{(n)} / n!$, we complete the proof of  lemma.

\begin{cor} Assume that $a_j \in \{0,1\}$ if $j\leq n$ and let the asterisk $*$ over a sum stand for the condition $a_j=1$. Then
\[
\gamma_n(k)= \theta^k\sum_{j_1\leq n}^{\ast}{1\over{j_1}}\cdots \sum_{j_k\leq
n}^{\ast}{{\mathbf{1}}\big\{j_1 +\cdots + j_k \leq  n\big\}\over{j_k}}
\psi_n\big(n-(j_1+\cdots+j_k)\big).
\]
\end{cor}

Afterwards, all error terms can depend on $\theta$. We will indicate dependence on other parameters if necessary. The symbol $\ll$ is used as an analog of $O(\cdot)$ and $a\asymp b$ means that $a\ll b$ and $b\ll a$.

\begin{lem} \label{lem-FMtr} If $a_j\in \Z_+\cap[0, m]$ for $j\leq n$, then
\begin{eqnarray}
 \gamma_n(k)
&=& \sum_{u=1}^k\theta^u\sum_{r_1+\cdots+r_u=k}{k-1\choose r_1-1}\cdots
      {k-r_1-\cdots-r_{u-1}-1 \choose r_u-1}\nonumber\\
      &&
      \quad \times\sum_{j_1+\cdots+j_u< n}{ a_{j_1(r_1)}\cdots
       a_{j_u(r_u)}\over{j_1\cdots j_u}}\bigg(1-{j_1 + \cdots + j_u\over n}\bigg)^{\theta-1} \nonumber\\
       &&\quad  +O\Big({1+\log^k n\over n^{1\wedge \theta}}\Big),
       \label{appr}
\end{eqnarray}
where $1\wedge \theta:=\min \{1,\theta\}$,  $n\geq 1$ and the constant in $O(\cdot)$ depends on  $m$ and $k$.
\end{lem}

{\it Proof}.  It suffices to deal with the case if $\theta\not=1$ and $n$ is sufficiently large. Set $\Delta_n(k)$ for the difference of $\gamma_n(k)$ in (\ref{gammank}) and the main term in its approximation (\ref{appr}). Using (\ref{gammank}) and the given bound of $a_j$, we have
\begin{eqnarray*}
   \Delta_n(k)
   &\ll&
  \sum_{u=1}^k C_u(k,m)\sum_{j_1,\dots,j_u< n}{{\mathbf 1}\{j_1+\cdots+j_u<n\}\over{j_1\cdots j_u}}\, {1\over n-(j_1 + \cdots + j_u)}\\
  &&\quad \times\Big(1-{j_1 + \cdots + j_u\over n}\Big)^{\theta-1}\\
   &&\quad + n^{1-\theta}\sum_{u=1}^k C_u(k,m)
\sum_{j_1,\dots,j_u\leq n}{{\mathbf 1}\{j_1+\cdots+j_u=n\}\over{j_1\cdots j_u}}.
      \end{eqnarray*}
Here
\[
              C_u(k,m):=\sum_{r_1+\cdots+r_u=k\atop 1\leq  r_i\leq m, i\leq u}{{k-1}\choose{r_1-1}}\cdots
      {{k-r_1-\cdots-r_{u-1}-1}\choose{r_u-1}}\ll 1
\]
if $1\leq u\leq k$. Using the latter, we see that a typical sum to be estimated is
  \begin{eqnarray}
 &&\sum_{j_1<n}{1\over{j_1}}\cdots \sum_{j_u<
n} {{\mathbf 1}\{j_1+\cdots+j_u<n\}\over
j_u\big(n-(j_1+\cdots+j_u)\big)} \Big(1-{j_1+\cdots+j_u\over n}\Big)^{\theta-1}\nonumber\\
&&\quad +
n^{1-\theta}\sum_{j_1\leq n}{1\over{j_1}}\cdots \sum_{j_u\leq
n}{{\mathbf 1}\{j_1+\cdots+j_u=n\}\over{j_u}}=:R_{nu}+r_{nu}
  \label{Rnk}
  \end{eqnarray}
  where $1\leq u\leq k$.
Now, in the sums of second remainder, at least one $j_i\geq n/u$, $1\leq i\leq u$. Hence
\begin{eqnarray*}
r_{nu}&\leq& {u^2\over n^{\theta}}\sum_{j_1\leq n}{1\over{j_1}}\cdots \sum_{j_{u-1}\leq
n}{{\mathbf 1}\{j_1+\cdots+j_{u-1}\leq n-n/u\}\over j_{u-1}}\\
&\leq& {u^2\over n^{\theta}}\bigg(\sum_{j\leq n}{1\over j}\bigg)^{u-1}\ll {\log^{u-1} n\over n^{\theta}}
\end{eqnarray*}
for every $1\leq u\leq k$.

For brevity, introduce  temporarily the notation $J=j_1+\cdots +j_u$ and $j=j_{u+1}$.
We will apply the mathematical induction for either of the sums in the splitting
\begin{eqnarray*}
     &&R_{n,u+1}\ll\sum_{j_1<n}{1\over{j_1}}\cdots \sum_{j_u<n}{{\mathbf 1}\{J< n\}\over j_u}
         \sum_{j\leq (n-J)/2}{1\over{j}}\, {1\over (n-J)-j}\\
         &&+\sum_{j_1<n}{1\over{j_1}}\cdots \sum_{j_u< n}{{\mathbf 1}\{J< n\}\over j_u}\sum_{(n-J)/2< j< n-J}{1\over{j}}\, {1\over (n-J)-j}\Big(1-{J+j\over n}\Big)^{\theta-1}\\
         &=&:R_{n,u+1}'+ R_{n,u+1}''.
     \end{eqnarray*}
  Now,
  \begin{eqnarray*}
     R_{n1}'+ R_{n1}''&=&\sum_{j\leq n/2}{1\over j}{1\over n-j}
     +\sum_{n/2< j< n}{1\over j}{1\over n-j}\Big(1-{j\over n}\Big)^{\theta-1}
          \\
     &\ll&  {\log n\over n}+ {1\over n^{\theta}}\sum_{n/2<j<n}(n-j)^{\theta-2}\ll { \log n\over n}+{1\over n^{1\wedge \theta}}
     \ll {\log n\over n^{1\wedge \theta}}.
     \end{eqnarray*}

     Assuming that $R_{nu}'\ll (\log^u n)/n$, we have
     \[
       R_{n,u+1}'\ll R_{nu}' \log n\ll (\log^{u+1} n)/n
       \]
in either of the cases $\theta<1$ or $\theta>1$.
     Further, if $\theta>1$, then $(1-(J+j)/n)^{\theta-1}\leq 1$. An easy estimation of the most inner sum now implies
     \[
 R_{n,u+1}''\ll R_{nu}' \log n\ll (\log^{u+1} n)/n.
       \]

     If  $\theta<1$, then
     \begin{eqnarray*}
   R_{n,u+1}''&\ll& \sum_{j_1< n}{1\over{j_1}}\cdots \sum_{j_u<n}{{\mathbf 1}\{J< n\}\over j_u(n-J)}\\
   &&\quad\times
     \sum_{(n-J)/2< j< n-J}\Big(1-{J+j\over n}\Big)^{\theta-1}\, {1\over (n-J)-j}\\
     &\ll&
     {1\over n^{\theta-1}}\sum_{j_1<n}{1\over{j_1}}\cdots \sum_{j_u<n}{{\mathbf 1}\{J< n\}\over j_k(n-J)}
     \sum_{1\leq s<n} s^{\theta-2}\\
     &\ll& {1\over n^{\theta-1}} {\log^u n\over n}=  {\log^u n\over n^\theta}
\end{eqnarray*}
since the last sum is bounded and the remaining iterated sum  has been estimated.

  Collecting all the estimates, we return to (\ref{Rnk}) and conclude that $R_{nu}+r_{nu}\ll(\log^u n) n^{-\theta\wedge 1}$ for sufficiently large $n$. Inserting this into expression of $\Delta_n(k)$, we furnish the proof of  lemma.

\s

In a similar way, we can follow after the growth of the factorial moments of $w(\sigma,J_n)$, i.e. that of $ \upsilon_n(l)$ as $l\to\infty$.

\begin{lem} \label{lem-upsilon} Let $J_n\subset \{1,\dots, n\}$ be arbitrary.
If $\theta\geq1$, then  $\upsilon_n(l)\leq \upsilon_n(1)^l$ for every $n, l\in \N$. If $\theta<1$, then
there exists a positive constant $C$ depending on $\theta$ only such that
  \[
                       \upsilon_n(l)\leq C^l  \big(\upsilon_n(1)+1\big)^l
  \]
  for every $l\in \N$.
  \end{lem}

  \textit{Proof}. The proof of the first claim is straightforward. In the case $\theta<1$, we apply the induction.
 Examine the most inner sum on the right-hand side  of the inequality
 \begin{equation}
 \upsilon_n(l+1)\leq \theta^l\sum_{j_1,\dots,j_l\in J_n}{{\mathbf 1}\{S<n\}\over j_1 \cdots j_l}\sum_{j\in J_n}{{\mathbf 1}\{j<n-S\}\over j}
 \Big(1 - {S + j\over n}\Big)^{\theta - 1},
 \label{recur-ups}
 \end{equation}
 where temporary $S:=j_1+\cdots+j_l$ and $j:=j_{l+1}$. The summands over $j\leq (n-S)/2$ contribute not more than
 \[
      2^{1-\theta}(1-S/n)^{\theta-1} \upsilon_n(1)
\]
and
\begin{eqnarray*}
&&\sum_{j\in J_n}{{\mathbf 1}\{(n-S)/2< j<n-S\}\over j}
 \Big(1 - {S + j\over n}\Big)^{\theta - 1}\\
 &\leq&
  {2\over n^{\theta-1}(n-S)}\sum_{j\in J_n}{\mathbf 1}\{(n-S)/2< j<n-S\}(n-S-j)^{\theta-1}\\
 &\leq&
 {2\over n^{\theta-1}(n-S)}\sum_{k<(n-S)/2} k^{\theta-1}\leq C_1 \Big(1-{S\over n}\Big)^{\theta-1}.
 \end{eqnarray*}
The last two estimates and (\ref{recur-ups}) yield
\[
   \upsilon_n(l+1)\leq  (2\vee C_1)\big(\upsilon_n(1)+1\big) \upsilon_n(l).
\]
 Consequently, the desired claim hold with $C=2\vee C_1:=\max\{2, C_1\}$.

\medskip

Let us introduce the concentration function
$$
Q_{n}(u) = \sup_{x \in \mathbf{R}} \nu_{n}(|h(\sigma)- x| < u), \quad u \geq 0, \; x\in\R,
$$
and
$$
D_{n}(u;\lambda) = \sum_{j \leq n} \frac{u^2 \wedge
(a_j-\lambda j)^2}{j},\qquad D_{n}(u) = \min_{\lambda \in
\mathbf{R}} D_{n}(u;\lambda).
$$

\begin{lem}
\label{lem-CF}
We have
\begin{equation}
Q_{n}(u) \ll uD_{n}(u)^{-1/2}
\label{CF}
\end{equation}
for every $\theta>0$.
\end{lem}

 \textit{Proof.} See \cite{TK-LMD09}.

The last lemma is used to obtain lower estimates of the further needed frequencies.
Let $J \subset\{j:j\leq n\}$ be an arbitrary
nonempty set, maybe, depending on $n$, and $\overline{J}=\{j:j\leq
n\}\setminus J$.

\begin{lem} \label{lem-JK} Let $\theta\geq1$,  $K>0$, and $J$ be such that
\begin{equation}
            \sum_{j\in J}{1\over j}\leq K<\infty.
\label{K}
\end{equation}
 Denote
 \[
 \mu_n(K)=\inf_{J}\, \nu_{n}\big(k_j(\sigma) = 0 \;\forall\,j \in J\big),
 \]
 where the infimum is taken over $J$ satisfying $(\ref{K})$. For a sufficiently large $n_0(K)$, there
  exists a positive constant $c(K)$, depending at most on $\theta$ and $K$, such that
$\mu_n(K)\geq c(K)$
if $n\geq n_0(K)$.

Moreover, for any $I\subset J\cap [1,n-n_0(K)]$ and
\[
\widetilde{S}_n:=
\bigcup_{j\in I} S_n^{j}:=\bigcup_{j\in I}\Big\{\sigma \in \S: k_j(\sigma)= 1,\;  k_i
(\sigma) = 0\;\forall i \in J\setminus \{j\}\Big\},
\]
we have that
\begin{equation}
\nu_{n}(\widetilde{S}_n)\geq c(K)\sum_{j\in I}{1\over j}\psi_n(n-j)\gg \sum_{j\in I}{1\over j}\Big(1-{j\over n}\Big)^{\theta-1}
  \label{nuSn}
\end{equation}
provided that $n\geq 2n_0(K)$.
\end{lem}

\proof The first claim is  Corollary 1.3 of Theorem 1.2 (see
\cite{TKEM-RIMS12}). The second claim is proved in \cite{TKEM-AUBud13}.

\medskip
The next observation supplies a possibility to apply the previous lemma of the sieve type.

\begin{lem} \label{prop-1} Assume that $h_n(\sigma)\in\Z$ and   $T_n(x)\Rightarrow  F(x)$. Then
   \begin{equation}
 \sum_{j \leq n} {{\mathbf 1}\{a_j\not=0\}\over j}\ll 1 \label{bound}
 \end{equation}
 provided that $\theta\geq1$ or $\theta<1$ and
   \begin{equation}
 \sum_{j \leq n} {{\mathbf 1}\{|a_j|\geq K\}\over j}\leq K_1 \label{bound1}
 \end{equation}
 for some positive constants $K$ and $K_1$. Now, the constant in $(\ref{bound})$ depends also on $F$, $K$, and $K_1$.
 \end{lem}

 {\it Proof}. Since the limit law has an atom,  we obtain a lower estimate of the  concentration function $Q_n(u)\geq c>0$ for every $u>0$ if $n$ is sufficiently large. Now applying Lemma \ref{lem-CF},
  we have $D_n(u,\lambda)\ll c^{-2}u^2$ for some $\lambda\in\R$. This, if $u\to 0$, yields the estimate
    \begin{equation}
 \sum_{j \leq n} {{\mathbf 1}\{a_j\not=\lambda j\}\over j}\ll 1. \label{bound2}
 \end{equation}
 Actually, $\lambda\in\Z$. Indeed, if $||\cdot||$ denotes the distance to the nearest integer, we have
  \[
  1\gg D_n(1,\lambda)\geq \sum_{j\leq n}{||\lambda j||^2\over j}
  \]
  and, further, $\lambda=:\tilde \lambda+\delta$, where $\tilde\lambda\in\Z$ and $\delta=O(n^{-1})$.
  Now the inequality $(x+y)^2\leq 2x^2+2y^2$, $x,y\in\R$, implies
 \[
     D_n(u,\tilde\lambda)\leq 2 D_n(u,\lambda)+2 \sum_{j\leq n} {u^2\wedge (\delta j)^2\over j}\ll u^2.
 \]
  Consequently, we may proceed with $\lambda\in\Z$.

  If $\theta\geq 1$, then, denoting $J:=\{j\leq n:\; a_j\not=\lambda j\}$ and
 \[
    h_n(\sigma)=\lambda \ell(\bar k(\sigma))+\sum_{j\in J} (a_j-\lambda j)k_j(\sigma)
  =:\lambda n+\tilde h_n(\sigma),
  \]
 by Lemma \ref{lem-JK}, for sufficiently large $n$,
\[
   \nu_n(h_n(\sigma)=\lambda n)= \nu_n(\tilde h_n(\sigma)=0)\geq \nu_{n}\big(k_j(\sigma) = 0 \;\forall\,j \in J\big)\geq c_1>0.
 \]
 Hence if $\lambda n\to\infty$ for some subsequence of $n\to\infty$, at least $c_1$ of the probability distribution mass of $h_n(\sigma)$ disappears at infinity. This contradicts to the assumption of theorem. Hence $\lambda\ll n^{-1}$ and, thus,  $\lambda=0$ if $n$ is sufficiently large.
 Now, the  estimate $D_n(1,0)\ll1$ contains (\ref{bound}).

 Assume  that $\theta<1$ is arbitrary and $|a_j|\leq K$  for the most part of $j\leq n$ in the sense of (\ref{bound1}). Now, manipulating with  the latter and
 the estimate (\ref{bound2}), we obtain the bound $\lambda\ll K/n$ which implies that $\lambda=0$ eventually.

The lemma is proved.

\smallskip

{\bf Remark}. If $a_j\in \Z$, $j\geq 1$, do not depend on $n$, then the additive function posses a limit distribution, i.e. $\nu_n(h(\sigma)<x)\Rightarrow F(x)$, if and only if the series
\[
                 \sum_{j\geq 1} {{\mathbf 1}\{a_j\not=0\}\over j}
\]
converges. The fact is well known \cite{ABT}, Theorem 8.25, since the three series in an analog of Kolmogorov's theorem for the integer-valued functions  reduce to this one. The last lemma gives a very short proof of the  necessity.

\section{Proof of Theorem 1}

Firstly we observe that, if the limit r.v is concentrated in $\Z_+$, we may confine ourselves to non-negative additive functions.  Set $a^+$ for the non-negative part of $a\in \R$. Let $h_n^{(+)}(\sigma)$ be the additive function defined as in (\ref{h}) via $a_j^+$, where $j\leq n$, and
\[
T_n^{(+)}(x):=\nu_n\big( h_n^{(+)}(\sigma)<x\big).
\]

\begin{prop} \label{prop-2} If $\theta\geq 1$, then the convergence  $T_n(x)\Rightarrow  F(x)$ is equivalent to $T_n^{(+)}(x)\Rightarrow F(x)$ together with condition $(\ref{ajneg})$.
    \end{prop}

\textit{Proof}. Assume that $T_n(x)\Rightarrow  F_Y(x)$. Then, by Lemma \ref{prop-1}, condition (\ref{bound}) holds.
 Set
\[
I:= \{j\leq n-n_0:\; a_{j} \leq -1\}\subset J:= \{j\leq n,\; a_{j} \not=0\},
\]
where $n\geq n_0$ and $n_0\in\N$, depending of $F$, is sufficiently large.
Define, as in Lemma \ref{lem-JK},
$$
S_{n}^j=\{\sigma \in \mathbf{S}_n:\; k_j(\sigma)=1, k_i(\sigma)=0 \quad \forall i\in J\setminus \{j\}\}
$$
and observe that $h_n(\sigma) = a_j \leq -1$ for all $\sigma \in S_{n}^j$ with $j\in I$. We have from Lemma \ref{lem-JK} and the above assumption that
$$
o(1) = \nu_n\big(h_n(\sigma) \leq -1\big) \geq \nu_n\bigg(\bigcup_{j \in I}S_{n}^j\bigg)\gg  \sum_{j \in I}\frac{1}{j}\Big(1-\frac{j}{n}\Big)^{\theta -1}.
$$
The sum can be extended also over $n-n_0\leq j\leq n$. This proves the necessity of (\ref{ajneg}).

Further, having (\ref{ajneg}), we claim that $T_n^{(+)}(x)\Rightarrow F_Y(x)$ is equivalent to
$T_n(x)\Rightarrow F_Y(x)$. Indeed,
\begin{eqnarray*}
\nu_n\big(h_n(\sigma)\not= h_n^{(+)}(\sigma)\big) &\leq&
\sum_{j\leq n \atop a_j\leq -1}\nu_n\big(k_j(\sigma) \geq 1\big)\\
&\leq&
\sum_{j\leq n \atop a_j\leq -1}\mathbf{E}_n k_j(\sigma) = \sum_{j\leq n \atop a_j\leq -1} \frac{\theta}{j}\psi_n(n-j)= o(1).
\end{eqnarray*}
We have used a particular case of formula (\ref{Wat}).

    Proposition \ref{prop-2} is proved.

    \medskip

\noindent\textit{Proof of Theorem} \ref{thm-but}. We have from Lemma \ref{prop-1} that the sum (\ref{bound}) is bounded by a constant $C_F$.  Further, by Proposition \ref{prop-2}, we see that condition (\ref{ajneg}) is satisfied and we may assume that $h_n(\sigma)\in \Z_+$.  In what follows the constants involved in estimates can depend on $F$.

  For an integer $m\geq 1$, we set  $a_j(m)=a_j$ if $a_j\leq m$ and $a_j(m)=m$ otherwise and introduce the truncated functions
\[
h_n(\sigma; m):=\sum_{j\leq n} a_{j}(m) k_j(\sigma).
\]
By Lemmas \ref{FM} and  \ref{lem-FMtr}, we have
\[
   \E_n h_n(\sigma,m)_{(l)}=\Upsilon_n(l,m)+{\rm o}(1).
\]
So, the purpose lays in proving that
\begin{eqnarray}
\lim_{m \rightarrow \infty} \limsup_{n \rightarrow \infty} \E_n h_n(\sigma,m)_{(l)} =
\lim_{m \rightarrow \infty} \liminf_{n \rightarrow \infty}  \E_n h_n(\sigma,m)_{(l)} = \E Y_{(l)}=:\Upsilon(l)
\label{purpose}
\end{eqnarray}
for each natural numbers $l \leq \alpha-1-\e$.

Set $J_n:=\{j\leq n:\; a_j\not=0\}$, then
\begin{eqnarray}
\mathbf{E}_{n}h_n(\sigma, m)_{(l)} &\leq&
 m^l \E_n w(\sigma,J_n)_{(l)}=m^l (\upsilon_n(l)+o(1))\nonumber\\
 &\leq& m^l(\upsilon_n(1)^l+o(1))\ll m^l
\label{limsup}
\end{eqnarray}
by virtue of  bound (\ref{bound}), where the hidden constant depends on $\theta, l$, and $F$.

We now split:
\begin{equation}
\mathbf{E}_{n}h_n(\sigma, m)_{(l)} =E_n(l,m)'+E_n''(l,m)+E_n'''(l,m),
\label{split}
\end{equation}
 where
$$
E_n'(l,m) = \sum_{b=l}^{m-1} b_{(l)} \nu_n \big(h_n(\sigma) = b\big),
$$
$$
E_n''(l,m) = \sum_{b=m}^{M} b_{(l)} \nu_n \big(h_n(\sigma,m) = b\big),
$$
$$
E_n'''(l,m) = \sum_{b>M} b_{(l)} \nu_n \big(h_n(\sigma,m) = b\big),
$$
and $M=M(m)>m$ is a natural number to be chosen later.

If  $l\leq \alpha$, then
\[
\lim_{n\to\infty}E_n'(l,m) = \sum_{b=l}^{m-1} b_{(l)}P (Y = b)=\Upsilon(l)-\sum_{b\geq m}b_{(l)}P (Y = b)
\]
for each fixed $m$ and, by virtue of  $\E Y^\alpha<\infty$,
 \[
 \lim_{m\to\infty}\sum_{b\geq m}b_{(l)}P (Y = b)\leq \lim_{m\to\infty} \sum_{b\geq m}b^l P (Y = b)=0.
 \]
 In other words,
 \begin{equation}
\lim_{m\to \infty} \lim_{n\to \infty} \E_n' (l,m)=\Upsilon(l)
 \label{E1}
 \end{equation}
 for each $l\leq \alpha$.

  Similarly, if $l\leq \alpha-1-\e$,
\begin{eqnarray}
E_n''(l,m) &\leq&  \sum_{b=m}^M b^l \nu_n\big(h_n(\sigma)\geq b\big) =
 \sum_{b=m}^M b^l P(Y\geq b)+ o_m(1)\nonumber\\
 &\leq&
  \E Y^{\alpha}\sum_{b=m}^{\infty} {1\over b^{1+\e}}+ o_m(1)=\rho_n(m).
 \label{E2}
\end{eqnarray}

Finally, we have from (\ref{limsup})
\begin{eqnarray*}
   E_n'''(l,m)&=&{1\over \theta^{(n)}}\sum_{\sigma\in\S}\theta^{w(\sigma)} {\mathbf 1}\{h_n(\sigma,m)>M\} h_n(\sigma,m)_{(l)} \cdot {h_n(\sigma,m)-l\over h_n(\sigma,m)-l}\\
   &\leq&
   {1\over M-l} \E_n h_n(\sigma, m)_{(l+1)}\ll {m^{l+1}\over M-l}\ll {1\over m}
   \end{eqnarray*}
for the choice $ M=m^{l+2}$ provided that $m>2l$.   Collecting (\ref{E1}), (\ref{E2}), and the last estimate, from the splitting (\ref{split}), we obtain claim (\ref{purpose}).

The theorem is proved.

\section{Proof of Theorem 2}

As we have seen in the proof of Proposition \ref{prop-2}, condition (\ref{ajneg}) allows us to deal with nonnegative functions only.
By the condition of theorem and Lemma \ref{lem-FMtr},
\[
\E_n h_n(\sigma,m)_{(l)}=\Upsilon_n(l,m)+o(1)=\Upsilon(l)+\rho_n(m).
 \]
Let $L\in \N$ be a fixed number, and examine the expansion of the characteristic function
\[
  \E_n {\re}^{it h_n(\sigma,m)}=\sum_{l=0}^L {\E_n h_n(\sigma,m)_{(l)}\over l!} (\re^{it}-1)^l+
  {O}\bigg(
  {\E_n h_n(\sigma,m)_{(l)}\over (L+1)!} |{\re}^{it}-1|^{L+1}\bigg),
\]
where $t\in \R$ and the constant in $O(\cdot)$ is absolute. We further have
\[
  \E_n {\re}^{it h_n(\sigma,m)}=\sum_{l=1}^L {\Upsilon(l)\over l!} ({\re}^{it}-1)^l+ {O}\bigg(
  {2^L \Upsilon(L+1)\over (L+1)!}\bigg) +\rho_m(n)
\]
uniformly in $t\in \R$. In other words,
$$
\lim_{m \to \infty}\limsup_{n\to\infty}\bigg|\mathbf{E}_n e^{ith(\sigma;m)} - \sum_{l=0}^L \frac{\Upsilon(l)}{l!}(e^{it} - 1)^l \bigg| \ll
\frac{2^L \Upsilon(L+1)}{(L+1)!}
$$
for every $L \geq 1$.

By virtue of the given conditions,
\begin{eqnarray*}
&&\limsup_{n \to \infty} \sum_{j \leq n \atop a_j > m} \frac{\theta}{j}\bigg(1 - \frac{j}{n}\bigg)^{\theta - 1} \\
&\leq& \frac{1}{m}
\lim_{r \to \infty}\limsup_{n \to  \infty} \sum_{j \leq n} \frac{\theta a_j(r)}{j}\bigg(1 - \frac{j}{n}\bigg)^{\theta - 1}
={\Upsilon(1)\over m}.
\end{eqnarray*}
Hence
\begin{eqnarray*}
\mathbf{E}_{n}| e^{it h(\sigma;m)} - e^{it h_n(\sigma)}| &\leq& \nu_n(h(\sigma) \neq h(\sigma;m))
\leq \sum_{j \leq n \atop a_j > m}\nu_n(k_j(\sigma) \geq 1)\\
 &\leq& \sum_{j \leq n \atop a_j > m}\E_n k_j(\sigma)= \sum_{j \leq n \atop a_j > m}\frac{\theta}{j}\psi_n(n-j)=\rho_m(n).
\end{eqnarray*}
The last two approximations, imply
$$
\lim_{m \to \infty}\limsup_{n \to  \infty}\bigg|\mathbf{E}_n e^{ith_n(\sigma)} - \sum_{l=0}^L \frac{\Upsilon(l)(e^{it} - 1)^l}{l!}\bigg| \ll \frac{2^L \Upsilon(L+1)}{(L+1)!}.
$$
It remains to take $L \to \infty$.

The theorem is proved.

\s

\textit{Proof of Corollary} \ref{cor-n-puse}. In the sufficiency part, it suffices to rewrite the factorial moments as follows:
\begin{eqnarray*}
\Upsilon_n(l,m)&=&
\theta\sum_{n/2<j< n}{a_{j(l)}(m)\over j}\Big(1-{j\over n}\Big)^{\theta-1}\\
&=&
\sum_{k=1}^m k_{(l)}\sum_{n/2<j<n} {\theta{\mathbf 1}\{a_j=k\}\over j}\Big(1-{j\over n}\Big)^{\theta-1}.
\end{eqnarray*}
Here we have to check if the inner sums can approach the Poisson probabilities. Since their sum over $k\geq0$ tends to $t_\theta(1)\geq 1-{\rm e}^{-\mu}$, this is possible. We may continue and get
\[
\Upsilon_n(l,m)
=\sum_{k=1}^m k_{(l)}\Big({\rm e }^{-\mu} {\mu^k\over k!}+o_k(1)\Big)
={\rm e }^{-\mu}\sum_{k=1}^m k_{(l)} {\mu^k\over k!}+o_m(1).
\]
Hence
\[
        \Upsilon_n(l,m)-\mu^l=\rho_n(m)
\]
as desired.

 To prove the necessity, we demonstrate another path. Recall that the  function $\varphi_n(z)=(\theta^{(n)}/n!)\E_n {\rm e}^{zh(\sigma)}$ satisfies (\ref{deriv}). If $a_j=0$ for $j\leq n/2$, then $\varphi_{n-j}(z)=\theta^{(n-j)}/(n-j)!$ and, consequently,  we obtain
  \begin{eqnarray*}
{\rm e}^{\mu(z-1)} +o(1)&=&  \E_n z^{h_n(\sigma)}=
1+  \theta\sum_{n/2<j\leq n}{z^{a_j}-1\over j}\psi_n(n-j)\\
&=&
1+  \theta\sum_{k\geq 1}(z^k-1)\sum_{n/2<j\leq n}{{\mathbf 1}\{a_j=k\}\over j}
\Big(1-{j\over n}\Big)^{\theta-1}+ o(1)
\end{eqnarray*}
uniformly in $z$ if $|z|=1$. Applying Cauchy's formula, we complete the proof of the corollary.

\s

\textit{Example.} Let $t_\theta(x)$ and $\mu$ be as in Corollary \ref{cor-n-puse}. Introduce the sequence $1/2=d_0<d_1<\cdots$ by
\[
        t_\theta(d_m)={\rm e}^{-\mu} \sum_{k=1}^m {\mu^k\over k!}, \quad m=0,1,\dots.
\]
and set $a_j=m$ if $nd_{m-1}<j\leq n d_m$  and $a_j=0$ otherwise. let $h_n(\sigma)$ be the completely additive function defined via these $a_j$. We claim that it posses the Poisson limit law with parameter $\mu$.

Firstly, we check that the function
is strictly increasing in $x$. We also observe that $t_\theta(1)<1$. Indeed, this is evident if $1\leq \theta\leq 1/\log 2$. Otherwise,
$   t_\theta(1)\leq t_{1/\log 2}(1)$ since it is decreasing in $\theta\geq 1/\log 2$.
The observed properties assure that the sequence  $d_m$ in the proposition is correctly defined. Moreover, approximating the sum by the Riemann integral, we have
\begin{eqnarray*}
&&\sum_{n/2<j<n} {\theta{\mathbf 1}\{a_j=m\}\over j}\Big(1-{j\over n}\Big)^{\theta-1}\\
&=&
\sum_{nd_{m-1}<j<nd_m} {\theta{\mathbf 1}\{a_j=m\}\over j}\Big(1-{j\over n}\Big)^{\theta-1}\\
&=&
t_\theta(d_m)-t_\theta(d_{m-1})+o(1)={\rm e}^{-\mu} {\mu^m\over m!}+o(1).
\end{eqnarray*}

The claim now follows from the last corollary.

\section{The cases with bounded $a_j$}

\textit{Proof of Theorem \ref{thm-aprez}}. Condition (\ref{aprez}) allows us to explore the case with $0\leq a_j\leq K$, $j\leq n$, only.

\textit{Sufficiency}. If $J_n:=\{j\leq n:\; a_j\not=0\}$, then
\[
           \E_n h_n(\sigma)_{(l)}\leq K^l \E_n w(\sigma,J_n)_{(l)}\leq K^l C^l\big(\upsilon_n(1)+1\big)^l \leq C_2^l
  \]
by Lemma \ref{lem-upsilon}  for every $l\geq1$. Further, it suffices to apply Theorem \ref{thm-pak}.

  \textit{Necessity}.  By Lemma \ref{prop-1}, we obtain from  the convergence $T_n(x)\Rightarrow F_Y(x)$ the bound (\ref{bound}) which, in its turn, yields $\Upsilon_n(1)\ll K$. Indeed, to check this, it suffices to observe that the summands over $n/2<j<n$ contribute only the bounded quantity. As we have seen,
in the sufficiency part,
\[
                 \sup_{n\geq1} \E_n h_n(\sigma)_{(l)}\ll K^l
\]
for every fixed $l\geq1$. Now, the weak convergence of distributions implies also the convergence of moments. Namely, we have $\Upsilon_n(l)=\E Y_{(l)}+o(1)$ where $l\geq 1$.

The theorem is proved.

\s

\textit{Proof of Corollary} \ref{cor-Poi}. Only \textit{Necessity} requires some argument. By Theorem \ref{Poisson}, convergence of distributions implies the relations $\upsilon_n(l)\to\mu^l$ where $l\geq1$.  Omitting nonnegative sums in the difference below, we obtain
\begin{eqnarray*}
o(1)&=&\upsilon_n(1)^l-\upsilon_n(l)\\
 &\geq&
\theta^l\sum_{j_1,\dots,j_l\leq n}^\ast{{\mathbf 1}\{j_1+\cdots+j_l>n\}\over{j_1\cdots j_l}}\Big(1 - {j_1\over n} \Big)^{\theta - 1}\cdots \Big(1 - {j_l\over n} \Big)^{\theta - 1}\\
&\geq&
\bigg(\theta^l\sum_{n/l<j\leq n}^\ast{1\over j}\Big(1 - {j\over n} \Big)^{\theta - 1}\bigg)^l
\end{eqnarray*}
for every $l\geq 1$. This yields the second  of conditions in (\ref{Poi-tr}). Using the latter and checking that the factor $(1-j/n)^{\theta-1}=1+o(1)$ uniformly in $j\leq r=o(n)$, we can rewrite the relation $\upsilon_n(1)=\mu+o(1)$ as is given in the first of relations in (\ref{Poi-tr}).

\s

\textit{Proof of Corollary} \ref{cor-finsup}. \textit{Sufficiency.} Since  (\ref{up-nL}) and (\ref{tr-nL}) imply the sufficient condition (\ref{upsilon}) in Corollary \ref{cor-Poi}, we have done.

 \textit{Necessity}. In the discussed case, the  $L$th factorial moment $\E_n h_n(\sigma)_{(L)}$  converges to zero. Hence the relevant formula  yields
 \begin{eqnarray*}
 o(1)&=&\E_n w(\sigma,J_n)_{(L)}\geq\theta^l\sum_{j_1,\dots j_L\leq n/L}^\ast{1\over{j_1\cdots j_L}}\psi_n\big(n-(j_1 + \cdots + j_L)\big)\\
 &\geq&
 \bigg(\theta\sum_{j\leq n/L}^\ast{1\over{j}}\psi_n(n-j)\bigg)^L
 \end{eqnarray*}
 for $\theta\leq 1$. This is equivalent to (\ref{tr-nL}). It also allows us to reduce the problem to the sequence of additive functions with $a_j=0$ if $j\leq n/L$. Then the necessary condition (\ref{upsilon}) reduces to  (\ref{up-nL}).

 The corollary is proved.

\section{Instances}

Let $\mathfrak L$ be the class of possible limit distributions for $w(\sigma,J_n)$, where $J_n\subset \{1,\dots, n\}$ is arbitrary, under the Ewens probability measure $\nu_n$. It would be desirable to find its description; we present some instances, however. We now apply Corollary \ref{cor-Y}.

\s

\textit{Bernoulli distribution} $\operatorname{Be}(p)$, where $p$ is the parameter $p\in (0,1)$. We claim that  $\operatorname{Be}(p)\in \mathfrak L$ if $p\leq t_\theta(1)$, where $t_\theta(x)$ is the previously defined function on $[1/2,1]$ and $\theta\geq 1$. The construction is based upon the factorial moments. For   $\operatorname{Be}(p)$, they are  $\upsilon(1) = p$ and $\upsilon(l) = 0$ if $ l \geq 2$. It suffices, therefore, to find $\alpha$ such that $t_\theta(\alpha)=p$ and to define $J_n=\{j\leq n:\; n/2< j\leq \alpha n\}$. By a simple approximation of the sum  by the integral, we verify  condition (\ref{upsilon}) and find that  $\upsilon_n(1)=p+o(1)$.
 \s

\textit{Binomial distribution} $\operatorname{Bi}(M, p)$, where $M\in\N$ and $p\in (0,1)$. Now, the factorial moments are equal to
$ M_{(l)}p^l$ if  $ l = 1,2,\ldots, M$ and to zero if $l = M+1, M+2, \dots$. For simplicity, we confine ourselves to a particular case of $\theta=1$ and $M=2$.
We claim that $\mathcal{B}(p,2)\in \mathcal L$ if
\[
0<p\leq (1/2)\log 3=0.405...
\]
The idea how to construct such an instance has been shown in \cite{JSGS-LMJ11}. Let us take two temporary parameters $0<\alpha\leq \log2$ and $0<\beta\leq \log(3/2)$. Define the sequence of  sets of natural numbers
\[
    J_n=\N\cap \Big((n/3, (n/3){\rm e}^\alpha]\cup
    ( 2n/3,  (2n/3){\rm e}^\beta]\Big).
\]
 The factorial moments of the additive function $w(\sigma,J_n)$ are equal to
\[
   \gamma_n(1)=  \sum_{j\in J_n} {1\over j}=\alpha+\beta+o(1)\leq \log 3+o(1),
    \]
\[
    \gamma_n(2)= \sum_{i,j\in J_n} {\mathbf 1}\{i+j\leq n\}{1\over i j}=\bigg(\sum_{n/3<j\leq
    (n/3){\rm e}^\alpha} {1\over j}\bigg)^2= \alpha^2+o(1),
\]
and $\gamma_n(l)=0$ if $l\geq 3$. To get the binomial distribution, we have to require that
\[
      2\alpha^2=(\alpha+\beta)^2.
\]
Hence
\[
               \alpha=(\sqrt2+1)\beta.
\]
Given $p\leq (1/2)\log 3$, we can choose $\beta$ and, consequently, $\alpha$ so that
    \[
                     2p=\alpha+\beta= (\sqrt2+2)\beta\leq \log 3.
    \]
Now, taking $\beta=(2-\sqrt2)p$ and $\alpha=p \sqrt2$, due to the condition  on $p$, we have finished.

\s

\textit{ A counterexample to Lugo's conjecture.} Examine $w(\sigma,J_n)$ where $ J_n=\{j:\; n/3<j\leq n/2\}$. A routine  approximation of sums by the Riemann integrals yields the following asymptotic formulas for the first two factorial moments:
\begin{eqnarray*}
    \upsilon_{n}(1)&=& \sum_{n/3<j\leq n/2} {\theta\over j}\Big(1-{j\over n}\Big)^{\theta-1}+o(1) \\
    &=&
     \theta\int_{1/3}^{1/2}(1-u)^{\theta-1}\frac{du}{u} + o(1)=:\lambda + o(1)
    \end{eqnarray*}
and
\begin{eqnarray*}
    \upsilon_{n}(2)
    &=&
    \theta^2\sum_{n/3<i,j\leq n/2} {1\over i j}\Big(1-{i+j\over n}\Big)^{\theta-1}+o(1)\\
    &=&
        \theta^2\int_{1/3}^{1/2}\int_{1/3}^{1/2}(1-u-v)^{\theta-1}{du dv\over uv} + o(1).
\end{eqnarray*}
Hence
\begin{eqnarray*}
&&\upsilon_{n}(1)^2-\upsilon_{n}(2)\\
\qquad&=&\theta^2\int_{1/3}^{1/2}\int_{1/3}^{1/2}\Big((1-u)^{\theta-1}(1-v)^{\theta-1}-(1-u-v)^{\theta-1}\Big){du dv \over uv}\\
&&\quad  + o(1)\geq c_2>0
 \end{eqnarray*}
 if $\theta>1$ and $n$ is sufficiently large. Observing also that $\upsilon_n(l)=0$ if $l\geq3$, by Corollary \ref{cor-Y}, we see that the function $w(\sigma,J_n)$ obeys a limit distribution but it is not the $(2,\lambda)$ quasi-Poisson.

\s

\textit{The laws outside $\mathfrak L$}. As we have stressed $\upsilon_n(l)\leq \upsilon_n(1)^l$ if $\theta\geq1$. Consequently, the inequality should be preserved by the laws in  $\mathfrak L$. Actually, this observation is due to J.\v Siaulys and G. Stepanauskas \cite{JSGS-SMS08}.
 The distributions like geometric with a parameter $p\in(0,1)$ or
a mixed Poisson distribution $F_Y(x)=\Pi(x;{\beta,\lambda,\tau})$ defined by the factorial  moments
\[
\mathbf{E}Y_{(l)}= \beta\lambda^l + (1 - \beta)\tau^l, \quad l= 1,2, \dots
\]
where $0<\beta<1$, $\lambda,\tau>0$, and $\lambda\not=\tau$, do not belong to  $\mathfrak L$ if $\theta\geq1$,

\bigskip

{\bf Concluding remark}. Most of the just presented results can be obtained for the generalized Ewens probability measure
\[
     \nu_{n,\Theta}(A):={1\over \Theta_n} \sum_{\sigma\in A}\theta_j^{w(\sigma)},
     \]
     where $0<c_3\leq \theta_j q^{-j}\leq C_3<\infty$ if $j\leq n$, $q\geq1$ is a fixed constant, and $\Theta_n$ is an appropriate normalization. In some cases, unfortunately, we have to assume that $c_3=1$. An analytic technique to deal with the value distribution of mappings defined on $\S$ with respect to $ \nu_{n,\Theta}$ was proposed by the second author \cite{EM-CPC02}. Later it was extended (see, for instance, \cite{EM-OJM09} and \cite{VZ-RJ11}) when $\theta_j q^{-j}$ satisfy some averaged conditions.  The asymptotic distributions under generalized Ewens measure of $h_n(\sigma)$ were treated in \cite{EM-OJM09}. The second author's paper \cite{EM-DMTCS12} provides an approximation in the total variation distance of the truncated  cycle vector by an appropriate vector with independent coordinates which is a basic tool for a probabilistic approach. The latter was already applied  in  \cite{BoEM-Ku12} to prove a functional limit theorem.

     The recent papers \cite{BaGr-JCT05}, \cite{BeUelVe-AAP11}, \cite{ErUel-Arx12}  discuss  cases with different behavior of $\theta_j$, e.g. $\theta_j= {\rm e}^{j^\gamma}$, $j\leq n$, where $0<\gamma<1$. Hopefully, the described method of factorial moments will be of  use in these cases.

\bigskip

Tatjana Bak\v sajeva

Department of Mathematics and Informatics

 Vilnius University

  Naugarduko str. 24, LT-03225 Vilnius

LITHUANIA

\medskip

 Eugenijus Manstavi\v cius

Institute of Mathematics and Informatics

 Vilnius University

Akademijos 4, LT-08663 Vilnius

LITHUANIA

e-mail: eugenijus.manstavicius{@}mif.vu.lt


\begin{thebibliography}{70}

\bibitem{BArKDa-Ar11} G.B. Arous, K. Dang, \textit{ On fluctuations of eigenvalues of random permutation matrices}, arXiv:1106.2108v1 (2011)

\bibitem{ABT} R. Arratia, A.D. Barbour, S. Tavar\'e, Logarithmic Combinatorial Structures: a Probabilistic Approach,  EMS Publishing House, Z\"urich (2003)

\bibitem{BaGr-JCT05} A.D. Barbour, B.L. Granovsky, \textit{Random combinatorial structures: the convergent case}, J. Combin. Theory A 109,  203--220 (2005)

\bibitem{BeUel-CMPh09} V. Betz, D. Ueltschi, \textit{ Spatial random permutations and infinite cycles}, Commun. Math. Phys. 285, 469--501 (2009)

\bibitem{BeUel-PTRF11} V. Betz, D. Ueltschi, \textit{Spatial random permutations with small cycle weights}, Probab. Theory. Rel. Fields 149, 191--222 (2011)

\bibitem{BeUel-EJP11} V. Betz, D. Ueltschi, \textit{Spatial random
permutations and Poisson-Dirichlet law of cycle lengths}, Electron. J. Probab. 16(41), P. 41, 1173--1192 (2011)

\bibitem{BeUelVe-AAP11} V. Betz, D. Ueltschi, Y. Velenik,
\textit{Random permutations with  cycle weights}, Ann. Appl. Probab. 21, 312--331 (2011)

\bibitem{BoEM-Ku12} K. Bogdanas, E. Manstavi\v cius, \textit{Stochastic processes on weakly logarithmic assemblies},
in: Anal. Probab. Methods Number Theory 5, Kubilius Memorial Volume,  A. Laurin\v cikas \textit{et al} (Eds), TEV, Vilnius, 69--80 (2012)

\bibitem{ErUel-Arx12} N.M. Ercolani, D. Ueltschi, \textit{Cycle structure of random permutations with cycle weights},  arXiv:1102.4796v2 (2012)

\bibitem{ET-ZWar65} P. Erd\"{o}s, P. Tur\'{a}n, \textit{On some problems of a statistical group theory} I, Z. Wahrsch. Verw. Gebiete 4, 175--186 (1965)

\bibitem{HKCS-SPA00} B. Hambly, P. Keevash, N. O\'{}Connell, D. Stark, \textit{The characteristic polynomial of a random permutation matrix},
Stoch. Process. Appl. 90, 335--346 (2000)

\bibitem{HNNZ-ar11} Ch. Hughes, J. Najnudel, A. Nikeghball, D. Zeindler, \textit{Random permutation matrices under the genereralized Ewens measure} arXiv:1109.5010v1 (2011)

\bibitem{TK-SMS07} T. Kargina, \textit{Additive
         functions on permutations and the Ewens Probability}, \v Siauliai Math. Semin. 10, 33--41  (2007)

\bibitem{TK-LMD09} T. Kargina, \textit{Asymptotic distributions of the number of restricted cycles in a random permutation}, Lietuvos matem. rink. Proc. LMS, 50, 420--425 (2009)

\bibitem{TKEM-RIMS12} T.Kargina, E.\ Manstavi\v cius, \textit{Multiplicative functions on $Z_+^n$
 and the Ewens Sampling Formula}, RIMS K$\hat{o}$ky$\hat{u}$roku Bessatsu, B34, 137--151  (2012)

\bibitem{TKEM-AUBud13} T. Kargina, E. Manstavi\v{c}ius, \textit{The law of large numbers with respect to Ewens probability}, Ann. Univ. Sci. Budapest., Sect. Comp.  39, 227--238 (2013)

\bibitem{ML-EJC09} M. Lugo, \textit{Profiles of permutations}, Electron. J. Combin. 16, R99 (2009)

\bibitem{ML-ArX09} M. Lugo,  \textit{The number of cycles of specified normalized length in permutations}, arXiv:0909.2909vI (2009)


\bibitem{EM-LMJ96} E. Manstavi\v cius, \textit{Additive and multiplicative functions on random permutations}, Lith. Math. J. 36, 400--408 (1996)

\bibitem{EM-CPC02} E. Manstavi\v cius, \textit{Mappings on decomposable combinatorial structures: analytic approach}, Combin. Probab. Computing 11,
61--78 (2002)

\bibitem{EM-APM02} E. Manstavi\v cius, \textit{Functional limit theorem for sequences of mappings on the symmetric group}, in: Anal. Probab. Methods Number Theory 3, A.\ Laurin\v cikas \textit{et al} (Eds), TEV, Vilnius, 175--187 (2002)

\bibitem{EM-LMJ05} E. Manstavi\v cius, \textit{The Poisson distribution for the linear statistics on random permutations}, Lith. Math. J. 45, 434--446 (2005)

\bibitem{EM-AMUO05} E. Manstavi\v cius, \textit{Discrete limit laws for additive functions on the symmetric group}, Acta Math. Univ. Ostraviensis 13, 47--55 (2005)

\bibitem{EM-RJ08} E. Manstavi\v cius, \textit{Asymptotic value distribution of additive function defined on the symmetric group}, Ramanujan J. 17, 259--280 (2008)

\bibitem{EM-OJM09} E. Manstavi\v cius, \textit{ An analytic method in probabilistic combinatorics}, Osaka J. Math. 46, 273--290 (2009)

\bibitem{EM-LMJ11} E. Manstavi\v cius, \textit{A limit theorem for additive functions defined on the symmetric group}, Lith. Math. J. 51, 211--237 (2011)

\bibitem{EM-DMTCS12} E. Manstavi\v cius, \textit{On total variation approximations for random
assemblies}, Discrete Math. Th. Comuter Sci. Proc., AofA-12, 97–108 (2012)

\bibitem{JS-LMJ96} J. \v{S}iaulys, \textit{Convergence to the Poisson law II. Unbounded strongly additive functions}, Lith. Math. J. 36, 393--404 (1996)

\bibitem{JS-LMJ98}  J. \v{S}iaulys, \textit{Convergence to the Poisson law. III. Method of moments}, Lith. Math. J. 38, 374--390 (1998)

\bibitem{JS-LMJ00} J. \v{S}iaulys, \textit{Factorial moments of distributions of additive functions}, Lith. Math. J. 40, 389--508 (2000)

\bibitem{JSGS-LMJ11}  J. \v{S}iaulys, G. Stepanauskas, \textit{Binomial limit law for additive prime indicators}, Lith. Math. J. 51, 562--572 (2011)

\bibitem{JSGS-SMS08} J. \v{S}iaulys, G. Stepanauskas, \textit{Some limit laws for strongly additive prime indicators}, \v{S}iauliai Math. Semin. 3,
 235--246 (2008)

\bibitem{KW-AP00} K.L. Wieand, \textit{Eigenvalue distributions of random permutation matrices}, Ann. Probab. 28, 1563--1587 (2000)

\bibitem{KW-JTP03} K.L. Wieand, \textit{Permutation matrices, wreath products, and the distribution of eigenvalues},  J. Theoret. Probab. 16, 599--623 (2003)

\bibitem{VZ-LMJ02} V. Zacharovas, \textit{The convergence rate to the normal law of a certain variable defined on random polynomials},
Lith. Math. J. 42, 88--107 (2002)

\bibitem{VZ-LMJ04} V. Zacharovas, \textit{Distribution of the logarithm of the order of a random permutation}, Lith. Math. J. 44, 296--327 (2004)

\bibitem{VZ-RJ11} V. Zacharovas, \textit{Voronoi summation formulae and multiplicative functions on permutations}, Ramanujan J. 24(3), 289--329 (2011)

\bibitem{DZ-EJP10} D. Zeindler, \textit{Permutation matrices and the moments of their characteristic polynomial}, Electron. J. Probab. 15, P 34, 1092--1118 (2010)

\end{thebibliography}
\end{document}